\documentclass{article} 

\title{The Wasserstein Proximal Gradient Algorithm}
\usepackage[nonatbib,final]{neurips_2020}

\author{
	Adil Salim\\
  Visual Computing Center\\ KAUST\\ adil.salim@kaust.edu.sa
  \And
	Anna Korba\\
	Gatsby Computational Neuroscience Unit\\ University College London\\a.korba@ucl.ac.uk
	\And
	Giulia Luise\\
	Computer Science Department\\ University College London\\ g.luise16@ucl.ac.uk
}

\usepackage{url}            
\usepackage{booktabs}       
\usepackage{amsfonts}       
\usepackage{nicefrac}       
\usepackage{microtype}      

\usepackage{subcaption}

\usepackage{graphicx}

 \usepackage[english]{babel} 
  \usepackage{times}		
\usepackage[T1]{fontenc}
\usepackage[utf8]{inputenc} 
\DeclareUnicodeCharacter{00A0}{~}
\usepackage{amssymb,amsmath,amscd,amsfonts,amsthm,bbm,mathrsfs,yhmath}
\usepackage[shortlabels]{enumitem}
\usepackage{algorithm}
\usepackage{algpseudocode}
\usepackage{hyperref}
\usepackage[sort,nocompress]{cite}
\usepackage[textwidth=2cm, textsize=footnotesize]{todonotes}  
\newcommand{\cP}{{\mathcal P}} 
 
\newcommand{\cG}{{\mathcal G}} 
\newcommand{\cE}{{\mathcal E}} 
\newcommand{\cH}{{\mathcal H}}

\theoremstyle{definition}
\newtheorem{theorem}{Theorem}
\newtheorem{lemma}[theorem]{Lemma}
\newtheorem{corollary}[theorem]{Corollary}
\newtheorem{proposition}[theorem]{Proposition}

\newtheorem{remark}{Remark}

\DeclareMathOperator{\prox}{prox}
\DeclareMathOperator{\jko}{JKO}

\newcommand{\eqdef}{:=} 

\newcommand{\cB}{{\mathcal B}} 
 
\newcommand{\E}{{{\mathbb E}}} 
\newcommand{\bR}{{\mathbb R}} 
 
\newcommand{\bE}{{\mathbb E}} 
\newcommand{\sX}{{\mathcal X}} 
 
\newcommand{\sZ}{{\mathsf Z}} 
\newcommand{\sx}{{\mathsf x}} 
\newcommand{\cL}{{{\mathcal L}}}

\newcommand{\ps}[1]{\langle #1 \rangle}

\DeclareMathOperator*{\argmin}{arg\,min}

\begin{document}

\maketitle

\begin{abstract}
    Wasserstein gradient flows are continuous time dynamics that define curves of steepest descent to minimize an objective function over the space of probability measures (\textit{i.e.}, the Wasserstein space). This objective is typically a divergence w.r.t. a fixed target distribution. In recent years, these continuous time dynamics have been used to study the convergence of machine learning algorithms aiming at approximating a probability distribution. However, the discrete-time behavior of these algorithms might differ from the continuous time dynamics. Besides, although discretized gradient flows have been proposed in the literature, little is known about their minimization power. In this work, we propose a Forward Backward (FB) discretization scheme that can tackle the case where the objective function is the sum of a smooth and a nonsmooth geodesically convex terms. Using techniques from convex optimization and optimal transport, we analyze the FB scheme as a minimization algorithm on the Wasserstein space. More precisely, we show under mild assumptions that the FB scheme has convergence guarantees similar to the proximal gradient algorithm in Euclidean spaces.

\end{abstract}

\section{Introduction}
The task of transporting an initial distribution $\mu_0$ to a target distribution $\mu_{\star}$ is common in machine learning. This task can be reformulated as the minimization of a cost functional defined over the set of probability distributions. Wasserstein gradient flows~\cite{ambrosio2008gradient} are suitable continuous time dynamics to minimize such cost functionals. These flows have found applications in various fields of machine learning such as reinforcement learning~\cite{richemond2017wasserstein,zhang2018policy}, sampling~\cite{wibisono2018sampling,durmus2019analysis,bernton2018langevin,cheng2017convergence} and neural networks optimization~\cite{chizat2018global,mei2019mean}. Indeed, Wasserstein gradient flows can be seen as the continuous limit of several discrete time machine learning algorithms.  The analysis of continuous time dynamics is often easier than the analysis of their discrete time counterparts. Therefore, many works focus solely on continuous time analyses of machine learning algorithms such as variants of gradient descent ~\cite{taghvaei2018accelerated,wang2019accelerated,chizat2018global,mei2019mean,duncan2019geometry,blanchet2018family}. Besides, although discretized Wasserstein gradient flows have been proposed in the literature~\cite{ambrosio2008gradient,jordan1998variational,maury2011handling,carlier2017splitting,bowles2015weak,wibisono2018sampling}, most of them have not been studied as minimization algorithms. 

In this paper, we focus on the resolution, by a discrete time algorithm, of a minimization problem defined on the set $\cP_2(\sX)$ of probability measures $\mu$ over $\sX = \bR^d$ such that $\int \|x\|^2 d\mu(x) < \infty$. More precisely, $\mu_\star$ is defined as a solution to
\begin{equation}\label{eq:composite}
	\min_{\mu \in \cP_2(\sX)} \cG(\mu) \eqdef \cE_F(\mu) + \cH(\mu), 
	\end{equation}
	where $\cE_F$ is a potential energy $\cE_F(\mu)=\int F(x)d\mu(x)$ tied to a smooth convex function $F:\bR^d \to \bR$, and $\cH$ is a nonsmooth term \textit{convex along the generalized geodesics} defined by the Wasserstein distance. The potential $\cE_F$ plays the role of a data fitting term whereas $\cH$ can be seen as a regularizer.

\paragraph{Motivation for studying the template problem~\eqref{eq:composite}.} Many free energy minimization problems can be cast as Problem~\eqref{eq:composite}, see \cite[Proposition 7.7]{agueh2011barycenters} or more generally~\cite[Section 9]{ambrosio2008gradient}. For instance, $\cH$ can be an internal energy~\cite[Example 9.3.6]{ambrosio2008gradient}. In particular, if $\cH$ is the negative entropy, then $\cG$ boils down to the Kullback Leibler divergence (up to an additive constant) w.r.t. the Gibbs measure $\mu_\star \propto \exp(-F)$. This remark has been used by several authors to study sampling tasks as minimizing Problem~\eqref{eq:composite}~\cite{wibisono2018sampling,durmus2019analysis,liu2017stein,bernton2018langevin,cheng2017convergence}. Another example where $\cH$ is an internal energy is the case where $\cH$ is a higher order entropy. In this case, $\mu_\star$ follows a Barenblatt profile. Moreover, in the context of optimization of infinitely wide two layers neural networks~\cite{mei2019mean,chizat2018global}, $\mu_{\star}$ denotes the optimal distribution over the parameters of the network. 
In this setting, $F(x)=-2\int k(x,y)d\mu_{\star}(y)$ is non convex and $\cH(\mu)=\int k(x,y)d\mu(x)d\mu(y)$ is an interaction energy~\cite[Example 9.3.4]{ambrosio2008gradient}, with $k$ depending on the activation functions of the network. 
Moreover, $\cG$ boils down to a Maximum Mean Discrepancy w.r.t. $\mu_{\star}$~\cite{arbel2019maximum} under a well-posedness condition. Alternatively, $F$ can be a regularizer of the distribution on the parameters of the network (see Appendix).
 
\textbf{Related works.}
Wasserstein gradient flows are continuous time dynamics that can minimize~\eqref{eq:composite}. Several time discretizations of such flows have been considered~\cite{santambrogio2017euclidean,ambrosio2008gradient,wibisono2018sampling}. However, these discretization schemes have  been mainly analyzed as numerical schemes to approximate the continuous gradient flow, rather than as optimization algorithms. 

Numerous optimization algorithms to solve~\eqref{eq:composite}, relying on different time-discretization schemes of the Wasserstein gradient flow of $\cG$, have been proposed previously. For instance,
~\cite{frogner2018approximate,zhang2018policy,santambrogio2017euclidean} rely on the implementation of the JKO (Jordan-Kinderlehrer-Otto) scheme~\cite{jordan1998variational}, which can be seen as a proximal scheme (\textit{i.e.}, \textit{backward discretization}) in the Wasserstein distance. In this case, each step of the algorithm relies on evaluating the JKO operator of $\cG$, exploiting efficient subroutines for this operator. 
When $\cG$ is smooth, some gradient descent algorithms over the $\cP_2(\sX)$ (\textit{i.e.}, \textit{forward discretizations}) have also been proposed~\cite{chewi2020gradient,liu2016stein,liu2017stein}. However, their analysis is notably challenging without further convexity assumptions on $\cG$. For example~\cite{chewi2020gradient}, which tackles the resolution of the Wasserstein barycenter problem (a nonconvex problem) via gradient descent, considers Gaussian input distributions. 
Another time-discretization of Wasserstein gradient flows can be found in the Langevin Monte Carlo (LMC) algorithm in the sampling literature. LMC can be seen as a splitting algorithm involving a gradient step for $\cE_F$ and an exact solving of the flow of $\cH$ (\textit{Forward-Flow discretization}~\cite{wibisono2018sampling}) to solve~\eqref{eq:composite} in the case where $\cH$ is the negative entropy. Indeed, the Wasserstein gradient flow of the negative entropy can be computed exactly in discrete time.\footnote{
 The addition of a gaussian noise corresponds to a process whose distribution is solution of the gradient flow of the relative entropy (the heat flow).
} Several works provide non asymptotic analyses of LMC in discrete time, see~\cite{durmus2019analysis,cheng2017convergence,bernton2018langevin,wibisono2019proximal,vempala2019rapid,ma2019there,pereyra2016proximal} among others. However, the convergence rate of the LMC algorithm does not match the convergence rate of the associated continuous time dynamics and LMC is a biased algorithm. 

A natural approach to minimize the sum of a smooth and a nonsmooth convex functions over a Hilbert space is to apply the proximal gradient algorithm~\cite{bau-com-livre11}, which implements a gradient step for the smooth term and a proximal step for the nonsmooth one. Indeed, many nonsmooth regularizers admit closed form proximity operators\footnote{see www.proximity-operator.net}. In this work we propose and analyze a splitting algorithm, the Forward-Backward (FB) scheme, to minimize the functional $\cG$ over $\cP_2(\sX)$. The FB scheme can be seen as an analogue to the proximal gradient algorithm in the Wasserstein space. This algorithm has already been studied (but no convergence rate has been established) in~\cite{wibisono2018sampling} in the particular case where $\cH$ is the negative entropy. Notably, ~\cite[Section F.1]{wibisono2018sampling} explains why the FB scheme should be unbiased, unlike Forward-Flow discretizations like LMC. 

More precisely, the proposed algorithm implements a forward (gradient step) for the smooth term $\cE_F$ and relies on the JKO operator for the nonsmooth term $\cH$~\cite{jordan1998variational} only, which plays the role of a proximity operator. This approach provides an algorithm with lower iteration complexity compared to the "full" JKO scheme applied to $\cE_F + \cH$.

\textbf{Contributions.} In summary, the Wasserstein gradient flows to minimize~\eqref{eq:composite} are well understood and modelize many machine learning algorithms. Discretized Wasserstein gradient flows have been proposed but little is known about their minimization power. In this work, we propose a natural optimization algorithm to solve Problem~\eqref{eq:composite}, which is a Forward-Backward discretization of the Wasserstein gradient flow. This algorithm is a generalization of a discretization scheme proposed in~\cite{wibisono2018sampling}. Our main contribution is to prove non asymptotic rates for the proposed scheme, under the assumptions that $F$ is smooth and convex and that $\cH$ is convex along the generalized geodesics defined by the Wasserstein distance.  We show that the obtained rates fortunately match the ones of the proximal gradient algorithm over Hilbert spaces. 

The remainder is organized as follows. In Section~\ref{sec:preliminaries} we provide some background knowledge in optimal transport and gradient flows. In Section~\ref{sec:algorithm}, we introduce the Forward-Backward Euler discretization scheme. We 
study the FB scheme as an optimization algorithm and present our main results, \textit{i.e.}, non-asymptotic rates for the resolution of~\eqref{eq:composite} in Section~\ref{sec:analysis}. In Section~\ref{sec:experiments} we illustrate the performance of our algorithm for a simple sampling task. 
The convergence proofs are postponed to the appendix.

\section{Preliminaries}
\label{sec:preliminaries}
In this section, we introduce the notations and recall fundamental definitions and properties on optimal transport and gradient flows that will be used throughout the paper. 

\subsection{Notations}

In the sequel, $\cP_2(\sX)$ is the space of probability measures $\mu$ on $\sX$ with finite second order moment. Denote $\cB(\sX)$ the Borelian $\sigma$-field over $\sX$. For any $\mu \in \cP_2(\sX)$, $L^2(\mu)$ is the space of functions $f : (\sX,\cB(\sX)) \to (\sX,\cB(\sX))$ such that $\int \|f\|^2 d\mu < \infty$. Note that the identity map $I$ is an element of $L^2(\mu)$. For any $\mu \in \cP_2(\sX)$, we denote by $\|\cdot\|_\mu$ and $\ps{\cdot,\cdot}_\mu$ respectively the norm and the inner product of the space $L^2(\mu)$. For any measures $\mu,\nu$, we write $\mu \ll \nu$ if $\mu$ is absolutely continuous with respect to $\nu$, and we denote $Leb$ the Lebesgue measure over $\sX$. The set of regular distributions of the Wasserstein space is denoted by $\cP_2^r(\sX)\eqdef\{\mu\in \cP_2(\sX), \; \mu \ll Leb\}$. If $f,g : \sX \to \sX$, the composition $f\circ g$ of $g$ by $f$ is sometimes denoted $f(g)$.

\subsection{Optimal transport}
\label{sec:OT}
For every measurable map $T$ defined on $(\sX,\cB(\sX))$ and for every $\mu \in \cP_2(\sX)$, we denote $T_{\#}\mu$ the pushforward measure of $\mu$ by $T$ characterized by the `transfer lemma', i.e.:
\begin{equation}
\int \phi(y)dT_{\#}\mu(y)=\int \phi(T(x))d\mu(x)\quad \text{ for any measurable and bounded function }\phi.
\end{equation}
Consider the $2$-Wasserstein distance defined for every $\mu,\nu \in \cP_2(\sX)$ by
\begin{equation}
\label{eq:wass}
W^2(\mu,\nu) \eqdef \inf_{\upsilon \in \Gamma(\mu,\nu)} \int \|x-y\|^2d\upsilon(x,y),
\end{equation}
where $\Gamma(\mu,\nu)$ is the set of couplings between $\mu$ and $\nu$~\cite{villani2008optimal}, \textit{i.e.} the set of nonnegative measures $\upsilon$ over $\sX \times \sX$ such that $P_{\#} \upsilon = \mu$ (resp. $Q_{\#} \upsilon = \nu$) where $P : (x,y) \mapsto x$ (resp. $Q : (x,y) \mapsto y$) is the projection onto the first (resp. second) component. 

We now recall the well-known Brenier theorem~\cite{brenier1991polar},\cite[Section 6.2.3]{ambrosio2008gradient}.
\begin{theorem}\label{th:brenier_th} Let $\mu \in \cP_2^r(\sX)$ and $\nu \in \cP_2(\sX)$. Then,
	\begin{enumerate}
		\item There exists a unique minimizer $\upsilon$ of~\eqref{eq:wass}. Besides, there exists a uniquely determined $\mu$-almost everywhere (a.e.) map $T_{\mu}^\nu : \sX \to \sX$ such that $ \upsilon = (I, T_{\mu}^\nu)_{\#}\mu$ where $(I,T_{\mu}^\nu) : (x,y) \mapsto (x,T_{\mu}^\nu(x))$. Finally, there exists a convex function $f : \sX \to \bR$ such that $T_{\mu}^\nu = \nabla f$ $\mu$-a.e.
		\item As a corollary, 
		\begin{equation}\label{eq:brenier_formula}
			W^2(\mu,\nu) = \int \|x-\nabla f(x)\|^2d\mu(x) = \inf_{T:\,\,T_{\#}\mu=\nu}\int \|x-T(x)\|^2d\mu(x).
		\end{equation}
		\item If $g : \sX \to \bR$ is convex, then $\nabla g$ is well defined $\mu$-a.e. and if $\nu = \nabla g_{\#} \mu$, then $T_{\mu}^\nu = \nabla g$ $\mu$-a.e.
		\item If $\nu \in \cP_2^r(\sX)$, then $T_{\nu}^\mu \circ T_{\mu}^\nu = I$ $\mu$-a.e. and $T_{\mu}^\nu \circ T_{\nu}^\mu = I$ $\nu$-a.e.
	\end{enumerate}	
\end{theorem}
Under the assumptions of Theorem~\ref{th:brenier_th}, the map $T_{\mu}^{\nu}$ is called the optimal transport (OT) map from $\mu$ to $\nu$. In this paper, as it is commonly the case in the literature, we may refer to the space of probability distributions $\cP_2(\sX)$ equipped with the $2$-Wasserstein distance as the Wasserstein space.

\subsection{Review of Gradient Flows and their discretizations}

\subsubsection{In an Euclidean space}

Assume that $\sX$ is a Euclidean space, consider a proper lower semi-continuous function $G : \sX \to (-\infty,+\infty]$ and denote by $D(G) = \{x \in \sX, \;G(x) < \infty\}$ its domain. We assume that $G$ is convex, \textit{i.e.}, for every $x,z \in \sX$ and for every $\varepsilon \in [0,1]$, we have:
\begin{equation}
\label{eq:def-cvx}
G(\varepsilon z + (1-\varepsilon) x) \leq \varepsilon G(z) + (1-\varepsilon)G(x).
\end{equation}
Given $x \in \sX$, recall that $y \in \sX$ is a subgradient of $G$ at $x$ if for every $z \in \sX$,
$$
G(x) + \ps{y,z - x} \leq G(z).
$$
The (possibly empty) set of subgradients of $G$ at $x$ is denoted by $\partial G(x)$, and the map $x \mapsto \partial G(x)$ is called the subdifferential. If $G$ is differentiable at $x$, then $\partial G (x)=\{\nabla G(x)\}$ where $\nabla G(x)$ is the gradient of $G$ at $x$. The subdifferential of the convex function $G$ allows to define the gradient flow of $G$: for every initial condition $\sx(0) = a$ such that $\partial G(a) \neq \emptyset$, there exists a unique absolutely continuous function $\sx : [0,+\infty) \to \sX$ that solves the differential inclusion~\cite[Th. 3.1]{bre-livre73}, \cite[Th. 2.7]{pey-sor-10}
\begin{equation}
\label{eq:ode}
\sx'(t) \in \partial G(\sx(t)).
\end{equation}
One can check that the gradient flow of $G$ is also characterized by the following system of Evolution Variational Inequalities (EVI) :
$$
\forall z \in D(G),\quad \frac{d}{dt} \|\sx(t)-z\|^2 \leq -2\left(G(\sx(t)) - G(z)\right).
$$
In contrast to \eqref{eq:ode}, the former characterization allows to define the gradient flow without using the notion of subdifferential, a property that can be practical in nonsmooth settings. Moreover, the non-asymptotic analysis of discretized gradient flows in the optimization literature often relies on discrete versions of the EVI.\\


The existence of Gradient Flows can be established as the limit of a proximal scheme~\cite[Th. 2.14]{pey-sor-10},~\cite[Th. 5.1]{bia-hac-sal-(sub)jca17} when the step-size $\gamma \to 0$. Defining the proximity operator of $G$ as:
\begin{equation}
\label{eq:prox}
\prox_{\gamma G}(x) \eqdef \argmin_{y \in \sX} G(y) + \frac{1}{2\gamma}\|y-x\|^2,
\end{equation}
the proximal scheme is written as
\begin{equation}\label{eq:prox_scheme}
x_{n+1} = \prox_{\gamma G}(x_n),
\end{equation}
which corresponds to the proximal point algorithm to minimize the function $G$, see~\cite{martinet1970breve}. The proximal scheme can be seen as a \textit{Backward Euler discretization} of the gradient flow. Indeed, writing the first order conditions of \eqref{eq:prox_scheme}, we have
$$
x_{n+1} \in x_n - \gamma \partial G(x_{n+1}),\quad \text
{ or equivalently }\quad
\frac{x_{n+1} - x_n}{\gamma} \in - \partial G(x_{n+1}).
$$
Hence, each iteration of the proximal scheme requires  
solving an equation which can be intractable in many cases. The \textit{Forward Euler scheme} is a more tractable integrator of the gradient flow of $G$, but is less stable and requires the differentiability of $G$. Under this assumption, this scheme is written
\begin{equation}
\label{eq:forward-euler-scheme}
\frac{x_{n+1} - x_n}{\gamma} = -\nabla G(x_{n})\quad \text{ or equivalently }\quad x_{n+1} = x_n - \gamma \nabla G(x_{n}),
\end{equation}
which corresponds to the well-known gradient descent algorithm to minimize the function $G$.
Consider now the case where the function $G$ can be decomposed as $G = F + H$, where $F$ is convex and smooth and $H$ is convex and nonsmooth. To integrate the gradient flow of $G = F+H$, another approach is to use the Forward and the Backward Euler schemes for the smooth term and nonsmooth term respectively~\cite{pey-sor-10}. This approach is also motivated by the fact that in many situations, the function $H$ is simple enough to implement its proximity operator $\prox_{\gamma H}$. If $G = F+H$, the Forward Backward Euler scheme is written as
\begin{equation}
\label{eq:forward-backward-euler-scheme}
\frac{x_{n+1} - x_n}{\gamma} \in -\nabla F(x_{n}) -\partial H(x_{n+1}).
\end{equation}
Recalling the definition of the proximity operator, this scheme can be rewritten as
\begin{equation}
\label{eq:prox-grad}
x_{n+1} = \prox_{\gamma H}(x_n - \gamma \nabla F(x_n)),
\end{equation}
which corresponds to the proximal gradient algorithm to minimize the composite function $G$.

\subsubsection{In the Wasserstein space}

Consider a proper lower semi-continuous functional $\cG : \cP_2(\sX) \to (-\infty,+\infty]$ and denote $D(\cG) = \{\mu \in \cP_2(\sX),\; \cG(\mu) < \infty\}$ its domain. We assume that $\cG$ is convex along generalized geodesics defined by the 2-Wasserstein distance~\cite[Chap. 9]{ambrosio2008gradient}, \textit{i.e.} for every $\mu, \pi \in \cP_2(\sX), \nu \in \cP_2^r(\sX)$ and for every $\varepsilon \in [0,1]$, 
\begin{equation}
\label{eq:def-general-geo-cvx}
\cG((\varepsilon T_{\nu}^{\pi} + (1-\varepsilon)T_{\nu}^{\mu})_{\#}\nu) \leq \varepsilon \cG(\pi) + (1-\varepsilon)\cG(\mu),
\end{equation}
where $T_{\nu}^{\pi}$ and $T_{\nu}^{\mu}$ are the optimal transport maps from $\nu$ to $\pi$ and from $\nu$ to $\mu$ respectively. Note that the curve $\varepsilon \mapsto (\varepsilon T_{\nu}^{\pi} + (1-\varepsilon)T_{\nu}^{\mu})_{\#}\nu$ interpolates between $\mu$ ($\varepsilon = 0$) and $\pi$ ($\varepsilon = 1$). Moreover, if $\nu = \pi$ or $\nu = \mu$, then this curve is a geodesic in the Wasserstein space. Given $\mu \in \cP_2(\sX)$, $\xi \in L^2(\mu)$ is a strong Fréchet subgradient of $\cG$ at $\mu$ \cite[Chap. 10]{ambrosio2008gradient} if for every $\phi \in L^2(\mu)$,
$$
\cG(\mu) + \varepsilon\ps{\xi,\phi}_{\mu} + o(\varepsilon) \leq \cG((I + \varepsilon \phi)_{\#} \mu).
$$
The (possibly empty) set of strong Fréchet subgradients of $\cG$ at $\mu$ is denoted $\partial \cG(\mu)$. The map $\mu \mapsto \partial \cG(\mu)$ is called the strong Fréchet subdifferential. 
Conveniently, the strong Fréchet subdifferential enables to define the gradient flow of the functional $\cG$ \cite[Chap. 11]{ambrosio2008gradient}.
However in the nonsmooth setting that will be considered in this paper, the characterization of gradient flows through EVI will be more practical. 
The gradient flow of $\cG$ is the solution of the following system of EVI \cite[Th. 11.1.4]{ambrosio2008gradient}:
$$
\forall \pi \in D(\cG),\quad \frac{d}{dt} W^2(\mu(t),\pi) \leq -2\left(\cG(\mu(t)) - \cG(\pi)\right).
$$
We shall perform a non-asymptotic analysis of a discretized gradient flow scheme to minimize the functional $\cG$. Our approach is to prove a discrete EVI for this scheme. \\



The existence of gradient flows can be established as the limit of a minimizing movement scheme~\cite[Th. 11.2.1]{ambrosio2008gradient},~\cite{jordan1998variational}. Defining the JKO operator of $\cG$ as: 
\begin{equation}
\label{eq:jko}
\jko_{\gamma \cG}(\mu) \eqdef \argmin_{\nu \in \cP_2(\sX)} \cG(\nu) + \frac{1}{2\gamma} W^2(\nu,\mu),
\end{equation}
the JKO scheme is written as
$$
\mu_{n+1} \in \jko_{\gamma \cG}(\mu_n).
$$
The JKO operator can be seen as a proximity operator by replacing the Wasserstein distance by the Euclidean distance. Moreover, the JKO scheme can be seen as a Backward Euler discretization of the gradient flow. Indeed, under some assumptions on the functional $\cG$, using~\cite[Lemma 10.1.2]{ambrosio2008gradient} we have
$$
\frac{T_{\mu_{n+1}}^{\mu_n}-I}{\gamma}\in \partial \cG(\mu_{n+1}).
$$
Using Brenier's theorem, since $T_{\mu_{n+1}}^{\mu_n}\circ T_{\mu_n}^{\mu_{n+1}}=I$ $\mu_n$-a.e., there exists a strong Fréchet subgradient of $\cG$ at $\mu_{n+1}$ denoted by $\nabla_W \cG(\mu_{n+1})$ such that
$$
\mu_{n+1} = \left(I - \gamma \nabla_W \cG(\mu_{n+1})\circ T_{\mu_n}^{\mu_{n+1}}\right)_{\#}\mu_n.
$$
Each iteration of the JKO scheme thus requires the minimization of a function which can be intractable in many cases. As previously, the Forward Euler scheme is more tractable and enjoys additionally a simpler geometrical interpretation. Assume $\partial \cG(\mu) = \{\nabla \cG(\mu)\}$ is a singleton for any $\mu \in D(\cG)$ (some examples are given \cite[Section 10.4]{ambrosio2008gradient}). The Forward Euler scheme for the gradient flow of $\cG$ is written:
\begin{equation}
	\label{eq:forward-euler-wass}
\mu_{n+1} = (I - \gamma \nabla \cG(\mu_n))_{\#}\mu_n,
\end{equation}
and corresponds to the iterations of the gradient descent algorithm over the Wasserstein space to minimize $\cG$. Although the Wasserstein space is not a Riemannian manifold, it can still be equipped with a Riemannian structure and interpretation~\cite{mccann2001polar,otto2001geometry}. In particular, the Forward Euler scheme can be seen as a Riemannian gradient descent where the exponential map at $\mu$ is the map $\phi \mapsto (I + \phi)_{\#}\mu$ defined on $L^2(\mu)$.

\section{The Forward Backward Euler scheme}
\label{sec:algorithm}



Recall that our goal is to minimize $\cE_F+\cH$,  
where $\cE_F(\mu)=\int F(x)d\mu(x)$ for any $\mu \in \cP_2(\sX)$, and $\cH$ is a nonsmooth functional. Throughout this paper, we assume the following on the potential function $F$: there exists $L,\lambda \ge 0$ such that
\begin{itemize}
	\item \textbf{A1.} $F$ is $L$-smooth \textit{i.e.} $F$ is differentiable and $\nabla F$ is $L$-Lipschitz continuous; for all $(x,y)\in \sX^2$:
	\begin{equation}
		F(y) \leq F(x) + \ps{\nabla F(x),y-x} + \frac{L}{2}\|x-y\|^2.
	\end{equation}		
	\item \textbf{A2.} $F$ is $\lambda$-strongly convex (we allow $\lambda = 0$); for all $(x,y)\in \sX^2$:
	\begin{equation}
		F(x) \leq F(y) - \ps{\nabla F(x),y-x} - \frac{\lambda}{2}\|x-y\|^2.
	\end{equation}		
\end{itemize}



Moreover, we assume the following on the function $\cH$.
\begin{itemize}
	\item \textbf{B1.} $\cH: \cP_2(\sX) \to (-\infty,+\infty]$ is proper and lower semi-continuous. Moreover, $D(\cH) \subset \cP_2^r(\sX)$.	
	\item \textbf{B2.} There exists $\gamma_0 >0$ such that  $\forall \gamma \in (0,\gamma_0)$, $\jko_{\gamma\cH}(\mu) \neq \emptyset$ for every $\mu \in \cP_2(\sX)$.
	\item \textbf{B3.} $\cH$ is convex along generalized geodesics.
\end{itemize}

Assumptions~\textbf{B1} and~\textbf{B2} are general technical assumptions, used
in~\cite{ambrosio2008gradient} (see~\cite[Eq. 10.1.1a, 10.1.1b]{ambrosio2008gradient}) that are satisfied in relevant cases. Then~\cite[Prop. 9.3.2, 9.3.5 and 9.3.9]{ambrosio2008gradient} gives broad examples of (potential, interaction and internal) energies 
satisfying \textbf{B3}, \textit{e.g.}, potential (resp. interaction) energies if the potential (resp. interaction) term is convex, and entropies (see Appendix and~\cite[Remark 9.3.10]{ambrosio2008gradient}). Moreover, $\cH$ is often nonsmooth, see \textit{e.g.}~\cite[Section 3.1.1]{wibisono2018sampling}. Therefore, we use a Forward Backward Euler scheme to integrate the gradient flow of $\cG$. 
Let $\gamma>0$ a step size. The proposed Forward Backward Euler scheme is written, for $n \ge 0$, as:
\begin{align}
\nu_{n+1} &= (I - \gamma \nabla F)_{\#} \mu_n \label{eq:forward}\\
\mu_{n+1} &\in \jko_{\gamma \cH}(\nu_{n+1}) \label{eq:backward}.
\end{align}
This scheme can be seen as a proximal gradient algorithm over the Wasserstein space to minimize the composite function $\cG = \cE_F + \cH$. 


\begin{remark} To our knowledge, the two cases where $\jko_{\gamma \cH}$ can be computed in closed form are the case where $\cH(\mu) = \int R d\mu$ and $R$ is a proximable function~\footnote{\textit{i.e.}, $R$ is a lower semi-continuous proper convex function whose proximity operator has a closed form}, see \cite{bowles2015weak}; and the case where $\cH$ is the negative entropy and $F$ is quadratic, see~\cite[Example 8]{wibisono2018sampling} and Section~\ref{sec:num}. However, there exist subroutines to compute JKOs of generic functionals, as well documented in the optimal transport and PDE literature (see the review of different strategies in \cite{santambrogio2017euclidean}). The implementation of the FB scheme relies on these subroutines, similarly to proximal splitting algorithms in optimization that rely on subroutines to compute the proximity operator. Moreover, as several proximity operators admit a close form\footnote{see www.proximity-operator.net} one can hope to be able to compute $\jko_{\gamma \cH}$ for simple functionals $\cH$. 
	
For instance, when $\cH$ is the negative entropy (defined by $\cH(\mu)= \int \log\left(\mu(x)\right)d\mu(x)$ if  $\mu \ll Leb$ with density $\mu$ and $\cH(\mu) \eqdef +\infty$ else), we conjecture that a technique similar to~\cite{genevay2016stochastic} can be applied. Again, for this particular functional, the JKO is known in closed form in the Gaussian case \cite[Example 8]{wibisono2018sampling}. Other works of interest for the proposed FB scheme have investigated efficient methods to implement the JKO of a generic functional with respect to the entropy-regularized Wasserstein distance \cite{peyre2015entropic}.
\end{remark}
In the next section, we study the non asymptotic properties of the FB scheme.

\section{Non asymptotic analysis}
\label{sec:analysis}
In this section, we provide rates for the Wasserstein proximal gradient algorithm. The main technical challenge is to handle the fact that $\cH$ is only convex along generalized geodesics (and not along any interpolating curve). We overcome this challenge by using the intuition that the Wasserstein proximal gradient algorithm is a discretization of the associated Wasserstein gradient flow. More precisely, the main step is to establish a discrete EVI to prove the rates. The proof differs from the one of the proximal gradient algorithm because the convexity inequality can be used with optimal transport maps only, see Lemma~\ref{lem:general-geo-cvx}. 

We consider a fixed step size $\gamma < 1/L$ and a probability distribution $\pi \in \cP_2(\sX)$. Our main result (Proposition~\ref{prop:evi}) combines several ingredients: the identification of the optimal transport maps between $\mu_n, \nu_{n+1}$ and $\mu_{n+1}$ (see Equations \eqref{eq:forward} and \eqref{eq:backward}), the proof of a generic lemma regarding generalized geodesic convexity (Lemma~\ref{lem:general-geo-cvx}) and a proof of a discrete EVI (Lemma~\ref{lemma:discreteEVI}).

\subsection{Identification of optimal transport maps}

Lemmas~\ref{lem:forward},\ref{lem:jko} identify the optimal transport maps from $\mu_n$ to $\nu_{n+1}$ and from $\nu_{n+1}$ to $\mu_{n+1}$ in the Forward Backward Euler scheme, as soon as the step size is sufficiently small. In particular, Lemma~\ref{lem:jko} is a consequence of~\cite[Lemma 10.1.2]{ambrosio2008gradient}.
\begin{lemma}
	\label{lem:forward}	
	Assume \textbf{A1}. Let $\mu \in \cP_2^r(\sX)$ and $\nu = (I - \gamma \nabla F)_{\#} \mu$. Then if $\gamma < 1/L$, the optimal transport map from $\mu$ to $\nu$ corresponds to
	$$T_{\mu}^{\nu} = I - \gamma \nabla F.$$
	Moreover, $\nu$ belongs to $\cP_2^r(\sX).$
\end{lemma}
\begin{lemma}
	\label{lem:jko}
	Assume \textbf{B1, B2}.
	Let $\nu \in \cP_2(\sX)$. If $\mu \in \jko_{\gamma \cH}(\nu)$, then $\mu \in D(\cH) \subset \cP_2^r(\sX)$ 
	and the optimal transport map $T_{\mu}^{\nu}$ from $\mu$ to $\nu$ satisfies
	$T_{\mu}^{\nu} \in I + \gamma \partial \cH(\mu)$.
	In other words, there exists a strong Fréchet subgradient at $\mu$ denoted by  $\nabla_W \cH(\mu)$ such that
	\begin{equation}
	\label{eq:resolvent2}
	T_{\mu}^{\nu} = I + \gamma \nabla_W \cH(\mu).
	\end{equation}
\end{lemma}
Using Lemmas~\ref{lem:forward},\ref{lem:jko}, if $\mu_0 \in \cP_2^r(\sX)$, then $\mu_n, \nu_n \in \cP_2^r(\sX)$ for every $n$ by induction.
This remark allows to consider optimal transport maps from $\mu_n$ and $\nu_n$ to any $\pi\in \cP_2(\sX)$.
The next lemma extends~\cite[10.1.1.B]{ambrosio2008gradient} to functionals $\cH$ convex along generalized geodesics.
\begin{lemma}
	\label{lem:general-geo-cvx}
	Assume \textbf{B1, B2, B3}.
	Let $\nu \in \cP_2^r(\sX)$, $\mu,\pi \in \cP_2(\sX)$ and $T_{\nu}^{\mu},T_{\nu}^{\pi}$ the optimal transport maps from $\nu$ to $\mu$ and from $\nu$ to $\pi$ respectively. If $\xi \in \partial \cH(\mu)$, then
	\begin{equation}
	\ps{\xi \circ T_{\nu}^{\mu},T_{\nu}^{\pi} - T_{\nu}^{\mu}}_\nu \leq \cH(\pi) - \cH(\mu).
	\end{equation}
\end{lemma}
Lemma~\ref{lem:general-geo-cvx} is natural, holds for any functional convex along generalized geodesics, and it is novel, to the best of our knowledge. The following results rely on this lemma. 

\subsection{A descent lemma}

Without using any convexity assumption on $F$, we first obtain a descent lemma. We denote $Y_{n+1} \eqdef T_{\nu_{n+1}}^{\mu_{n+1}}$ the optimal transport map between $\nu_{n+1}$ and $\mu_{n+1}$ in the Forward Backward Euler scheme \eqref{eq:forward}, \eqref{eq:backward}, and $X_{n+1} \eqdef Y_{n+1} \circ (I - \gamma \nabla F)$. Note that $X_{n+1}$ is a pushforward from $\mu_n$ to $\mu_{n+1}$.
\begin{theorem}[Descent] Assume $\mu_0 \in \cP_2^r(\sX)$, $\gamma<1/L$ and \textbf{A1, B1, B2, B3}. Then for $n\geq 0$, there exists a strong Fréchet subgradient at $\mu_{n+1}$ denoted by $\nabla_W \cH(\mu_{n+1})$ such that:
	\label{th:descent}
	\begin{equation*}
	\cG(\mu_{n+1}) \leq \cG(\mu_{n}) -\gamma\left(1- \frac{L\gamma}{2}\right)\| \nabla F + \nabla_W \cH(\mu_{n+1})(X_{n+1})\|_{\mu_n}^2,
	\end{equation*}
	where  we use the notation $\nabla_W \cH(\mu_{n+1})(X_{n+1})$ to denote $\nabla_W \cH(\mu_{n+1})\circ X_{n+1}$.
\end{theorem}
Hence, the sequence $(\cG(\mu_n))_n$ is decreasing as soon as the step-size is small enough. 

\subsection{Discrete EVI}
To prove a discrete EVI and obtain convergence rates, we need the additional convexity assumption \textbf{A2} on the potential function $F$. We firstly prove the two following lemmas.
\begin{lemma}\label{lemma:discreteEVI}
	\label{lem:evi-prox} Assume $\mu_0 \in \cP_2^r(\sX)$, $\gamma<1/L$, and \textbf{B1, B2, B3}. Then for $n\geq 0$ and $\pi \in \cP_2(\sX)$, there exists a strong Fréchet subgradient at $\mu_{n+1}$ denoted $\nabla_W \cH(\mu_{n+1})$ such that:
	$$
	W^2(\mu_{n+1},\pi) \leq W^2(\nu_{n+1},\pi) - 2\gamma\left(\cH(\mu_{n+1}) - \cH(\pi)\right) - \gamma^2 \|\nabla_W \cH(\mu_{n+1})\|_{\mu_{n+1}}^2.
	$$
\end{lemma}

\begin{lemma} Assume $\mu_0 \in \cP_2^r(\sX)$,  $\gamma\leq1/L$, and \textbf{A1, A2} with $\lambda \geq 0$. Then for $n \ge 0$, and $\pi \in \cP_2(\sX)$
	\label{lem:evi-grad}
	$$
	W^2(\nu_{n+1},\pi) \leq (1-\gamma\lambda)W^2(\mu_{n},\pi) - 2\gamma\left(\cE_F(\mu_n) - \cE_F(\pi)\right) + \gamma^2\|\nabla F\|_{\mu_n}^2.
	$$
\end{lemma}
 We can now provide a discrete EVI for the functional $\cG = \cE_F + \cH$.
\begin{proposition}[discrete EVI] 
	\label{prop:evi}
	Assume $\mu_0 \in \cP_2^r(\sX)$, $\gamma<1/L$, and \textbf{A1}--\textbf{B3} with $\lambda \geq 0$. Then for $n\geq 0$ and $\pi \in \cP_2(\sX)$, there exists a strong Fréchet subgradient at $\mu_{n+1}$ denoted by $\nabla_W \cH(\mu_{n+1})$ such that the Forward Backward Euler scheme verifies:
	\begin{equation}
		\label{eq:evi-discrete}
		W^2(\mu_{n+1},\pi) \leq (1-\gamma \lambda)W^2(\mu_{n},\pi) -2\gamma \left(\cG(\mu_{n+1})- \cG(\pi)\right).
	\end{equation}
\end{proposition}
\subsection{Convergence rates}
When the potential function $F$ is convex, we easily get rates from the discrete EVI inequality provided above. Theorem~\ref{th:evi} is a direct consequence of Proposition~\ref{prop:evi} by taking $\pi = \mu_\star$, and its corollaries provide rates depending on the strong convexity parameter of $F$. 
\begin{theorem}
	\label{th:evi} Assume $\mu_0 \in \cP_2^r(\sX)$, $\gamma<1/L,$ and \textbf{A1}--\textbf{B3} with $\lambda \geq 0$. Then for every $n\geq 0$,
	$$
	W^2(\mu_{n+1},\mu_\star) \leq (1-\gamma\lambda)W^2(\mu_{n},\mu_\star) -2\gamma( \cG(\mu_{n+1}) - \cG(\mu_\star)).
	$$
\end{theorem}
\begin{corollary}[Convex case rate] Under the assumptions of Theorem~\ref{th:evi}, for $n \ge 0$:
	\label{cor:cvx}
	$$
	\cG(\mu_n)-\cG(\mu_\star) \leq \frac{W^2(\mu_0,\mu_\star)}{2\gamma n}.
	$$
\end{corollary}
\begin{corollary}[Strongly convex case rate] Under the assumptions of Theorem~\ref{th:evi}, if $\lambda >0$, then for $n \ge 0$:
	\label{cor:scvx}
	$$
	W^2(\mu_{n},\mu_\star) \leq (1-\gamma\lambda)^n W^2(\mu_{0},\mu_\star).
	$$
\end{corollary}
Hence, as soon as $F$ is convex, we get sublinear rates in terms of the objective function $\cG$, while when $F$ is $\lambda$-strongly convex with $\lambda>0$, we get linear rates in the squared Wasserstein distance for the iterates of the Forward Backward Euler scheme.
	These rates match those of the proximal gradient algorithm in Hilbert space in the convex and strongly convex cases. 
	\cite{nesterov2018lectures}. Moreover, these rates are discrete time analogues of the continuous time rates obtained in~\cite[Th. 11.2.1]{ambrosio2008gradient} for the gradient flow of $\cG$. In the particular case where $\cH$ is the negative entropy, we can compare the convergence rates of the FB scheme to those of LMC. Although, the complexity of one iteration of the FB scheme is higher, the convergence rates of the FB scheme are faster, see \textit{e.g.}~\cite{durmus2019analysis}.

\section{Numerical experiments}\label{sec:experiments}
\label{sec:num}
\begin{figure}[ht!]
	\[
	  \begin{array}{cc}
	 \includegraphics[width=.4\linewidth]{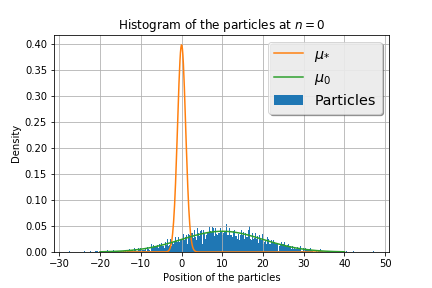} &
	  \includegraphics[width=.4\linewidth]{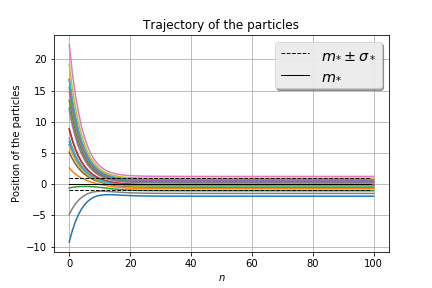}\\
		\includegraphics[width=.4\linewidth]{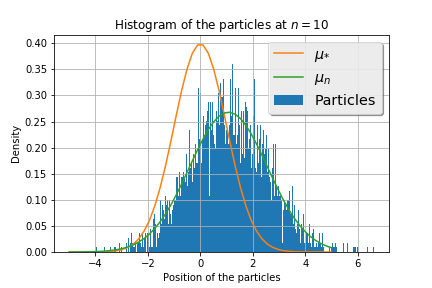} &
	  \includegraphics[width=.4\linewidth]{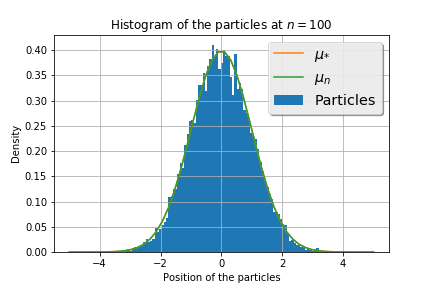}
	  \end{array}
	\]
		\caption{The particle implementation of the FB scheme illustrate the convergence of $\mu_n$ to $\mu_\star$.}
		\label{fig:particles}
	\end{figure}

We provide numerical experiments with a ground truth target distribution $\mu_{\star}$ to illustrate the dynamical behavior of the FB scheme, similarly to~\cite[Section 4.1]{taghvaei2018accelerated}. We consider $F(x) = 0.5|x|^2$, and $\cH$ the negative entropy. In this case, $\cG(\mu)$ is (up to an additive constant) the Kullback-Leibler divergence of $\mu$ w.r.t. the standard Gaussian distribution $\mu_\star$. We denote by $m_\star$ the mean and $\sigma_\star$ the variance, and fix
 $m_\star = 0$ and $\sigma_\star^2 = 1.0$. 
 We use the closed-form particle implementation of the FB scheme~\cite[Section G.1]{wibisono2018sampling}.
 This allows to show the dynamical behavior of the FB scheme when $\gamma = 0.1$, and $\mu_0$ is Gaussian with $m_0=10$ and $\sigma_0=100$, in Figure~\ref{fig:particles}. Note that $\lambda = 1.0$. 


 More precisely, the position of a set of particles initially distributed according to $\mu_0$ is updated iteratively. Using~\cite[Example 8]{wibisono2018sampling}, if $\mu_\star$ and $\mu_0$ are Gaussian, then $\mu_n$ is Gaussian for every $n$. Moreover, the mean (resp. the covariance matrix) of $\mu_{n+1}$ can be expressed as a function of the mean (resp. the covariance matrix) of $\mu_n$. Finally, the update of the position of the particles can be computed using these expressions, see the update formulas in~\cite[Example 8]{wibisono2018sampling}. These expressions also allow to compute $W^2(\mu_n,\mu_\star)$: since $\mu_n$ and $\mu_\star$ are Gaussian, $W^2(\mu_n,\mu_\star)$ can be computed in closed form and is also known as Bures-Wasserstein distance.

 The empirical distribution of the particles, represented by histograms, approximates $\mu_n$. We see that $\mu_n$ matches $\mu_\star$ after few iterations on Figure~\ref{fig:particles}. We also illustrate the particles, from their initial position (distributed according to $\mu_0$) to their final position (distributed according to $\mu_{100}$) in a high probability region of the target distribution $\mu_\star$. This shows that $\mu_{100}$ is close to $\mu_\star$. The linear convergence of $\mu_n$ to $\mu_\star$ in the squared Wasserstein distance is illustrated in a multidimensional case in Figure 2 (see Appendix).



\section{Conclusion}

We proposed an unified analysis of a Forward-Backward discretization of Wasserstein gradient flows in the convex case. We showed that the Forward-Backward discretization has convergence rates that are similar to the ones of the proximal gradient algorithm in Euclidean spaces.

Note that the implementation of the JKO operators is independent from the analysis of the FB scheme. However, an important problem raised by our work is whether we can find efficient implementations of the JKO operators of some specific functionals $\cH$ relevant in machine learning applications. 


\section{Acknowledgement}
Adil Salim thanks Pascal Bianchi for suggesting to use the JKO operator for optimization purposes. 

\section{Broader impact}
\label{sec:broader_impact}
The results that we showed, together with efficient implementations of some \textit{specific} JKOs could be very impactful for many machine learning tasks.

\bibliographystyle{plain}

\begin{thebibliography}{10}

\bibitem{agueh2011barycenters}
Martial Agueh and Guillaume Carlier.
\newblock Barycenters in the {W}asserstein space.
\newblock {\em SIAM Journal on Mathematical Analysis}, 43(2):904--924, 2011.

\bibitem{ambrosio2008gradient}
Luigi Ambrosio, Nicola Gigli, and Giuseppe Savar{\'e}.
\newblock {\em Gradient flows: in metric spaces and in the space of probability
  measures}.
\newblock Springer Science \& Business Media, 2008.

\bibitem{arbel2019maximum}
Michael Arbel, Anna Korba, Adil Salim, and Arthur Gretton.
\newblock Maximum mean discrepancy gradient flow.
\newblock In {\em Advances in Neural Information Processing Systems}, pages
  6481--6491, 2019.

\bibitem{bau-com-livre11}
Heinz~H Bauschke, Patrick~L Combettes, et~al.
\newblock {\em Convex analysis and monotone operator theory in Hilbert spaces},
  volume 408.
\newblock Springer, 2011.

\bibitem{bernton2018langevin}
Espen Bernton.
\newblock {L}angevin {M}onte {C}arlo and {JKO} splitting.
\newblock In {\em Conference on Learning Theory}, pages 1777--1798, 2018.

\bibitem{bia-hac-sal-(sub)jca17}
Pascal Bianchi, Walid Hachem, and Adil Salim.
\newblock A constant step {F}orward-{B}ackward algorithm involving random
  maximal monotone operators.
\newblock {\em Journal of Convex Analysis}, 2019.

\bibitem{blanchet2018family}
Adrien Blanchet and J{\'e}r{\^o}me Bolte.
\newblock A family of functional inequalities: {\L}ojasiewicz inequalities and
  displacement convex functions.
\newblock {\em Journal of Functional Analysis}, 275(7):1650--1673, 2018.

\bibitem{bowles2015weak}
Malcolm Bowles and Martial Agueh.
\newblock Weak solutions to a fractional {F}okker--{P}lanck equation via
  splitting and {W}asserstein gradient flow.
\newblock {\em Applied Mathematics Letters}, 42:30--35, 2015.

\bibitem{brenier1991polar}
Yann Brenier.
\newblock Polar factorization and monotone rearrangement of vector-valued
  functions.
\newblock {\em Communications on Pure and Applied Mathematics}, 44(4):375--417,
  1991.

\bibitem{bre-livre73}
Ha{\"i}m Br{\'e}zis.
\newblock {\em {Op\'erateurs maximaux monotones et semi-groupes de contractions
  dans les espaces de Hilbert}}.
\newblock North-Holland mathematics studies. Elsevier Science, Burlington, MA,
  1973.

\bibitem{carlier2017splitting}
Guillaume Carlier and Maxime Laborde.
\newblock A splitting method for nonlinear diffusions with nonlocal,
  nonpotential drifts.
\newblock {\em Nonlinear Analysis: Theory, Methods \& Applications}, 150:1--18,
  2017.

\bibitem{cheng2017convergence}
Xiang Cheng and Peter Bartlett.
\newblock Convergence of {L}angevin {MCMC} in {KL}-divergence.
\newblock In {\em Algorithmic Learning Theory}, pages 186--211, 2018.

\bibitem{chewi2020gradient}
Sinho Chewi, Tyler Maunu, Philippe Rigollet, and Austin~J Stromme.
\newblock Gradient descent algorithms for {B}ures-{W}asserstein barycenters.
\newblock In {\em Conference on Learning Theory}, pages 1276--1304, 2020.

\bibitem{chizat2018global}
Lenaic Chizat and Francis Bach.
\newblock On the global convergence of gradient descent for over-parameterized
  models using optimal transport.
\newblock In {\em Advances in Neural Information Processing Systems}, pages
  3036--3046, 2018.

\bibitem{duncan2019geometry}
Andrew Duncan, Nikolas N{\"u}sken, and Lukasz Szpruch.
\newblock On the geometry of stein variational gradient descent.
\newblock {\em arXiv preprint arXiv:1912.00894}, 2019.

\bibitem{durmus2019analysis}
Alain Durmus, Szymon Majewski, and Blazej Miasojedow.
\newblock Analysis of {L}angevin {M}onte {C}arlo via convex optimization.
\newblock {\em Journal of Machine Learning Research}, 20(73):1--46, 2019.

\bibitem{frogner2018approximate}
Charlie Frogner and Tomaso Poggio.
\newblock Approximate inference with {W}asserstein gradient flows.
\newblock In {\em International Conference on Artificial Intelligence and
  Statistics}, pages 2581--2590, 2020.

\bibitem{genevay2016stochastic}
Aude Genevay, Marco Cuturi, Gabriel Peyr{\'e}, and Francis Bach.
\newblock Stochastic optimization for large-scale optimal transport.
\newblock In {\em Advances in Neural Information Processing Systems}, pages
  3440--3448, 2016.

\bibitem{jordan1998variational}
Richard Jordan, David Kinderlehrer, and Felix Otto.
\newblock The variational formulation of the {F}okker--{P}lanck equation.
\newblock {\em SIAM Journal on Mathematical Analysis}, 29(1):1--17, 1998.

\bibitem{liu2017stein}
Qiang Liu.
\newblock Stein variational gradient descent as gradient flow.
\newblock In {\em Advances in Neural Information Processing Systems}, pages
  3115--3123, 2017.

\bibitem{liu2016stein}
Qiang Liu and Dilin Wang.
\newblock Stein variational gradient descent: A general purpose bayesian
  inference algorithm.
\newblock In {\em Advances in Neural Information Processing Systems}, pages
  2378--2386, 2016.

\bibitem{ma2019there}
Yi-An Ma, Niladri Chatterji, Xiang Cheng, Nicolas Flammarion, Peter Bartlett,
  and Michael~I Jordan.
\newblock Is there an analog of {N}esterov acceleration for {MCMC}?
\newblock {\em arXiv preprint arXiv:1902.00996}, 2019.

\bibitem{martinet1970breve}
Bernard Martinet.
\newblock Br{\`e}ve communication. r{\'e}gularisation d'in{\'e}quations
  variationnelles par approximations successives.
\newblock {\em Revue fran{\c{c}}aise d'informatique et de recherche
  op{\'e}rationnelle. S{\'e}rie rouge}, 4(R3):154--158, 1970.

\bibitem{maury2011handling}
Bertrand Maury, Aude Roudneff-Chupin, Filippo Santambrogio, and Juliette Venel.
\newblock Handling congestion in crowd motion modeling.
\newblock {\em Networks and Heterogeneous Media}, pages 485--519, 2011.

\bibitem{mccann2001polar}
Robert~J McCann.
\newblock Polar factorization of maps on riemannian manifolds.
\newblock {\em Geometric \& Functional Analysis GAFA}, 11(3):589--608, 2001.

\bibitem{mei2019mean}
Song Mei, Theodor Misiakiewicz, and Andrea Montanari.
\newblock Mean-field theory of two-layers neural networks: dimension-free
  bounds and kernel limit.
\newblock In {\em Conference on Learning Theory}, pages 1–--77, 2019.

\bibitem{nesterov2018lectures}
Yurii Nesterov.
\newblock {\em Lectures on convex optimization}, volume 137.
\newblock Springer, 2018.

\bibitem{otto2001geometry}
Felix Otto.
\newblock The geometry of dissipative evolution equations: the porous medium
  equation.
\newblock {\em Communications in Partial Differential Equations},
  26(1-2):101--174, 2001.

\bibitem{pereyra2016proximal}
Marcelo Pereyra.
\newblock Proximal {M}arkov chain {M}onte {C}arlo algorithms.
\newblock {\em Statistics and Computing}, 26(4):745--760, 2016.

\bibitem{pey-sor-10}
Juan Peypouquet and Sylvain Sorin.
\newblock Evolution equations for maximal monotone operators: asymptotic
  analysis in continuous and discrete time.
\newblock {\em Journal of Convex Analysis}, pages 1113--1163, 2010.

\bibitem{peyre2015entropic}
Gabriel Peyr{\'e}.
\newblock Entropic approximation of {W}asserstein gradient flows.
\newblock {\em SIAM Journal on Imaging Sciences}, 8(4):2323--2351, 2015.

\bibitem{richemond2017wasserstein}
Pierre~H Richemond and Brendan Maginnis.
\newblock On {W}asserstein reinforcement learning and the {F}okker-{P}lanck
  equation.
\newblock {\em arXiv preprint arXiv:1712.07185}, 2017.

\bibitem{santambrogio2017euclidean}
Filippo Santambrogio.
\newblock $\{$Euclidean, metric, and Wasserstein$\}$ gradient flows: an
  overview.
\newblock {\em Bulletin of Mathematical Sciences}, 7(1):87--154, 2017.

\bibitem{taghvaei2018accelerated}
Amirhossein Taghvaei and Prashant~G Mehta.
\newblock Accelerated gradient flow for probability distributions.
\newblock In {\em International Conference on Machine Learning}, 2019.

\bibitem{vempala2019rapid}
Santosh Vempala and Andre Wibisono.
\newblock Rapid convergence of the unadjusted {L}angevin algorithm:
  Isoperimetry suffices.
\newblock In {\em Advances in Neural Information Processing Systems}, pages
  8092--8104, 2019.

\bibitem{villani2008optimal}
C{\'e}dric Villani.
\newblock {\em Optimal transport: old and new}, volume 338.
\newblock Springer Science \& Business Media, 2008.

\bibitem{wang2019accelerated}
Yifei Wang and Wuchen Li.
\newblock Accelerated information gradient flow.
\newblock {\em arXiv preprint arXiv:1909.02102}, 2019.

\bibitem{wibisono2018sampling}
Andre Wibisono.
\newblock Sampling as optimization in the space of measures: The {L}angevin
  dynamics as a composite optimization problem.
\newblock In {\em Conference on Learning Theory}, page 2093–3027, 2018.

\bibitem{wibisono2019proximal}
Andre Wibisono.
\newblock Proximal {L}angevin algorithm: Rapid convergence under isoperimetry.
\newblock {\em arXiv preprint arXiv:1911.01469}, 2019.

\bibitem{zhang2018policy}
Ruiyi Zhang, Changyou Chen, Chunyuan Li, and Lawrence Carin.
\newblock Policy optimization as {W}asserstein gradient flows.
\newblock In {\em International Conference on Machine Learning}, pages
  5737--5746, 2018.

\end{thebibliography}
\newcommand{\noop}[1]{} \def\cprime{$'$} \def\cdprime{$''$} \def\cprime{$'$}


\newpage
\appendix



\clearpage

\section{Further numerical experiment}
We consider solving Problem~\eqref{eq:composite}, in a case where the ground truth $\mu_\star$ is a high dimensional Gaussian distribution~\cite{chewi2020gradient}.

We shall illustrate the linear convergence of the FB scheme this problem.

We consider a multidimensional extension of the simulation of Section~\ref{sec:num}. More precisely, $F(x) = 0.5\|x\|^2$ and, using the notations of Section~\ref{sec:num}, the inital distribution is $\mu_0^{\otimes d}$ and the target distribution is $\mu_\star^{\otimes d}$, where $\otimes$ denotes the product of measures.

The mean and covariance matrix of $\mu_n$ can be computed in closed form using~\cite[Example 8]{wibisono2018sampling}, therefore we can also compute $W^2(\mu_n,\mu_\star)$.

\begin{figure}[ht!]
	\[
	  \begin{array}{c}
     \includegraphics[width=\linewidth]{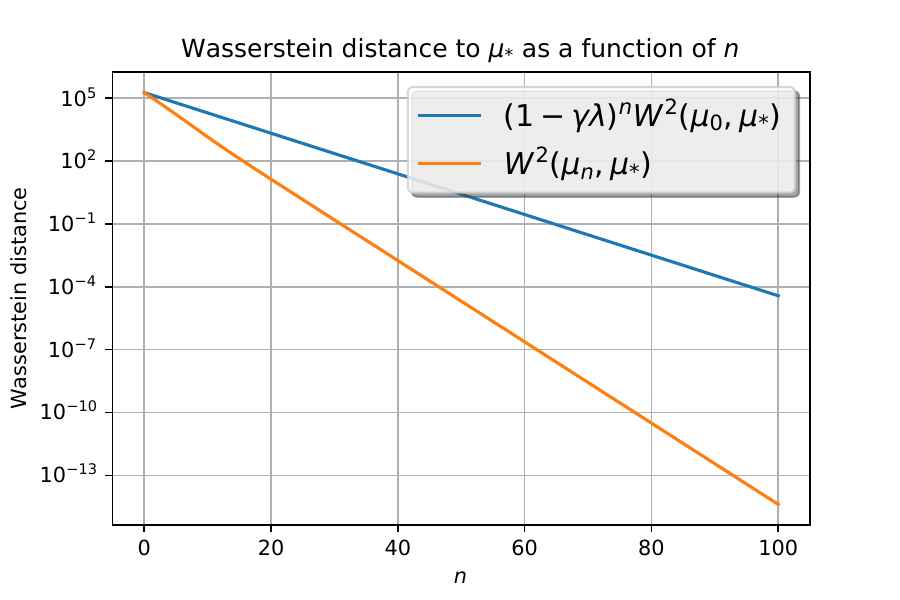} 
	  \end{array}
	\]
\caption{Linear convergence of $W^2(\mu_n,\mu_\star)$ to $0$ in dimension $d = 1000$.}
		\label{fig:lin}
	\end{figure}

The simulation can no longer be represented with histograms and particles, however, we represent the linear convergence of $\mu_n$ to $\mu_\star$ predicted by Corollary~\ref{cor:scvx} in Figure~\ref{fig:lin}.

We observe that $\mu_n$ converges linearly to $\mu_\star$ as predicted. Moreover, the observed linear rate is better than the one predicted by our theory. This suggests that our results can be improved in some particular cases.

\section{Examples of functionals $\cH$}

	\subsection{Internal energies}
	
	An internal energy is written 
	\begin{equation}
		\cH(\mu) = \int h(\mu(x))dx,
	\end{equation}
	if $\mu \ll Leb$ with density $\mu(x)$, and $\cH(\mu) = +\infty$ else. Internal energies satisfying Assumptions \textbf{B1}, \textbf{B2}, \textbf{B3} include the negative entropy, which corresponds to $h(u) = u \log(u)$, and the higher order entropies, which correspond to $h(u) = u^v / v$ where $v > 1$, see~\cite[Remark 9.3.10]{ambrosio2008gradient}. Moreover, in the case of the negative entropy, there is only one minimizer $\mu_\star$ of Problem~\eqref{eq:composite}, and its density w.r.t. $Leb$ can be written $\exp(C - F(x))$, where $C \in \bR$. In the case of the higher order entropies, $\mu_\star(x) = \max(0,C - F(x))^{\frac{1}{v-1}}$.
\subsection{[Nonconvex] Neural network optimization}

	Consider a infinite-width one hidden layer neural network (1HLNN)~\cite{arbel2019maximum,chizat2018global,mei2019mean} and let $n$ be the number of neurons. There are several ways to cast the optimization of a 1HLNN as Problem~\eqref{eq:composite}. 

	For any input $x$, the output of the 1HLNN can be written as $f(x)=\frac{1}{n}\sum_{i=1}^{n}\phi(x,z_i)$,  where $\phi(x,z_i)=w_i \psi(x,\theta_i)$ with $w_i\in \mathbb{R}$ the weight of the $i$-th neuron and $\theta_i\in \mathbb{R}
^d$ parametrizes the activation function $\psi$ (e.g. a sigmoid).  Given data $(x,y)\sim p$, the optimization of this neural network can be written as the minimization of the risk: $\argmin_{z_1, \dots, z_n \in \sZ }\mathbb{E}_{(x,y)\sim p}[\ell\left(y,\frac{1}{n}\sum_{i=1}^{n}\phi(x,z_i)\right)]$, where $\sZ=\bR \times \bR^d$ and $\ell$ is some loss function. When $n\rightarrow \infty$, this becomes an optimization problem over $\cP_2(\sZ)$: $\argmin_{\mu \in \cP(\sZ)}\bE_{(x,y)\sim p}[\ell\left(y,\int\phi(x,z)d\mu(z)\right)]$. Denote $R(h) = \bE_{(x,y)\sim p}[\ell\left(y,h(x)\right)]$.

	In~\cite{chizat2018global}, the optimization of a 1HLNN is cast as Problem~\eqref{eq:composite}, where $F$ is a regularizer on the parameters of the network (\textit{e.g.} $F(x) = \|x\|^2$) and $\cH(\mu) = R\left(\int \phi d\mu\right)$ is the risk. 
	
	Alternatively, consider the case where $\ell(y,x) = \frac{1}{2}\|y - x\|^2$ as in~\cite{arbel2019maximum}. Expanding the risk function leads to the objective $\cG(\mu)=\int F(z) d\mu(z)+ \frac{1}{2} \int K(z,z')d\mu(z)d\mu(z')$, where
$F(z)=-\E_{(x,y)\sim p}[y\phi(x,z)]$ is a potential and $K(z,z')=\E_{x\sim p_x}[\phi(x,z)\phi(x,z')]$ is an interaction term.

Note however that optimizing a 1HLNN is in general only a geodesically semiconvex  problem~\cite{arbel2019maximum,chizat2018global} (which is weaker than~\eqref{eq:def-general-geo-cvx}, see~\cite[Definition 9.1.1]{ambrosio2008gradient}), and hence is not strictly covered by our theory.

\section{Proof of Lemma~\ref{lem:forward}}

The map $I-\gamma \nabla F$ is a pushforward from $\mu$ to $\nu$. Moreover, denoting $u : x \mapsto \frac12\|x\|^2 - \gamma F(x)$, $\nabla u = I - \gamma \nabla F$.

By elementary algebra, for any $(x,y)\in \sX^2$ we have:
\begin{equation}
\label{eq:str-cvx-square}
\frac{1}{2}\|x\|^2 = \frac{1}{2}\|y\|^2 - \ps{x,y-x} - \frac{1}{2}\|x-y\|^2,
\end{equation}
and from the smoothness of $F$,
\begin{equation}
\label{eq:smooth-F}
F(y) \leq F(x) + \ps{\nabla F(x),y-x} + \frac{L}{2}\|x-y\|^2.
\end{equation}
Therefore, combining \eqref{eq:str-cvx-square} and \eqref{eq:smooth-F} multiplied by $\gamma$ gives:
\begin{equation}
\label{eq:grad-cvx}
\frac{1}{2}\|x\|^2 - \gamma F(x) \leq \frac{1}{2}\|y\|^2 - \gamma F(y) -\ps{x - \gamma \nabla F(x),y-x} - \frac{1}{2}(1 - L\gamma)\|x-y\|^2.
\end{equation}
In other words,
\begin{equation}
    \label{eq:grad-cvx-u}
    u(x) \leq u(y) -\ps{\nabla u(x),y-x} - \frac{1}{2}(1 - L\gamma)\|x-y\|^2.
    \end{equation}
Therefore, if $\gamma \leq 1/L$, then $u$ is convex and $\nabla u = T_{\mu}^\nu$ using Brenier's theorem. Moreover, if $\gamma < 1/L$ then $u$ is $(1 - L\gamma)$-strongly convex. In consequence, 
$$(1 - L\gamma)\|x-y\|^2 \leq \ps{x-y,\nabla u(x) - \nabla u(y)}.$$
Therefore, $\nabla u$ is injective. Furthermore, using the strong convexity of $u$ and~\cite[Lemma 5.5.3]{ambrosio2008gradient} (see also~\cite[Th. 6.2.3, Th. 6.2.7]{ambrosio2008gradient}), $\nu \in \cP_2^r(\sX)$. \\

\section{Proof of Lemma~\ref{lem:jko}}

Let $\mu \in \jko_{\gamma \cH}(\nu)$. Since $D(\cH) \subset \cP_2^r(\sX)$,~\cite[Lemma 10.1.2]{ambrosio2008gradient} implies $\mu \in D(\cH)$ and $\frac{1}{\gamma}(T_{\mu}^{\nu}-I)\in \partial \cH(\nu)$ is a strong subgradient of $\cH$ at $\nu$.

\section{Proof of Lemma~\ref{lem:general-geo-cvx}}
Since $\xi \in \partial \cH(\mu)$, for every $\phi \in L^2(\mu)$,
$$
\cH((I+\varepsilon\phi)_{\#}\mu) \geq \cH(\mu) + \varepsilon \ps{\xi,\phi}_\mu + o(\varepsilon).
$$
Applying the last inequality to $\phi = T_{\nu}^{\pi} \circ T_{\mu}^{\nu} - I$ and using the transfer lemma ($\mu = {T_{\nu}^{\mu}}_{\#}\nu$) we have
$$
\cH((T_{\nu}^{\mu}+\varepsilon(T_{\nu}^{\pi} - T_{\nu}^{\mu}))_{\#}\nu) = \cH((I+\varepsilon\phi)_{\#}\mu), $$
$$\ps{\xi,\phi}_\mu = \ps{\xi(T_{\nu}^{\mu}),T_{\nu}^{\pi} - T_{\nu}^{\mu}}_\nu,
$$
and 
\begin{equation}
\label{eq:frechet}
\frac{\cH((T_{\nu}^{\mu}+\varepsilon(T_{\nu}^{\pi} - T_{\nu}^{\mu}))_{\#}\nu) - \cH(\mu)}{\varepsilon} \geq \ps{\xi(T_{\nu}^{\mu}),T_{\nu}^{\pi} - T_{\nu}^{\mu}}_\nu + o(1).
\end{equation}
Using the generalized geodesic convexity of $\cH$, 
$$
\cH((T_{\nu}^{\mu}+\varepsilon(T_{\nu}^{\pi} - T_{\nu}^{\mu}))_{\#}\nu) \leq \varepsilon \cH(\pi) + (1-\varepsilon) \cH(\mu).
$$
Plugging the last inequality into~\eqref{eq:frechet},
\begin{equation}
\cH(\pi) - \cH(\mu) \geq \ps{\xi(T_{\nu}^{\mu}),T_{\nu}^{\pi} - T_{\nu}^{\mu}}_\nu + o(1).
\end{equation}
We get the conclusion by letting $\varepsilon \to 0$.

\section{Proof of Theorem~\ref{th:descent}}

Denote $Y_{n+1} \eqdef T_{\nu_{n+1}}^{\mu_{n+1}}$ the optimal transport map between $\nu_{n+1}$ and $\mu_{n+1}$ and $\nabla_W \cH(\mu_{n+1})$ the strong Fréchet subgradient of $\cH$ evaluated at $\mu_{n+1}$ defined by Lemma~\ref{lem:jko}: $T_{\mu_{n+1}}^{\nu_{n+1}} = I + \gamma \nabla_W \cH(\mu_{n+1})$. 
Since $\mu_{n+1}, \nu_{n+1} \in \cP_2^r(\sX)$, $(I + \gamma \nabla_W \cH(\mu_{n+1}))\circ Y_{n+1} = I$ using Brenier's theorem. We thus have $\nu_{n+1}$-a.e.:
\begin{equation}
\label{eq:Y_n}
Y_{n+1} = I - \gamma \nabla_W \cH(\mu_{n+1})(Y_{n+1}).
\end{equation}

We firstly bound the $\cH$ term. By taking $\mu=\mu_{n+1}, \pi=\mu_n$ and $\nu=\nu_{n+1}$ in Lemma~\ref{lem:general-geo-cvx}, we have:
\begin{equation}
\label{eq:H-cvx-descent}
\cH(\mu_{n+1}) \leq \cH(\mu_{n}) - \ps{\nabla_W \cH(\mu_{n+1})(T_{\nu_{n+1}}^{\mu_{n+1}}), T_{\nu_{n+1}}^{\mu_{n}} - T_{\nu_{n+1}}^{\mu_{n+1}}}_{\nu_{n+1}}.
\end{equation}
We now identify $T_{\nu_{n+1}}^{\mu_{n}}$ and $T_{\nu_{n+1}}^{\mu_{n+1}}$. Recall that $Y_{n+1} = T_{\nu_{n+1}}^{\mu_{n+1}}$. Moreover, using Brenier's theorem and Lemma~\ref{lem:forward}, $\nu_{n+1} \in \cP_2^r(\sX)$ and $T_{\nu_{n+1}}^{\mu_{n}} = (I - \gamma \nabla F)^{-1}$. Therefore, 
\begin{equation*}
\cH(\mu_{n+1}) \leq \cH(\mu_{n}) - \ps{\nabla_W \cH(\mu_{n+1})(Y_{n+1}), (I-\gamma \nabla F)^{-1} - Y_{n+1}}_{\nu_{n+1}}.
\end{equation*}
Using the transfer lemma, with $Y_{n+1}=X_{n+1}\circ (I-\gamma \nabla F)^{-1}$, the last inequality is equivalent to
\begin{equation}
\label{eq:descent-H}
\cH(\mu_{n+1}) \leq \cH(\mu_{n}) - \ps{\nabla_W \cH(\mu_{n+1})(X_{n+1}), I - X_{n+1}}_{\mu_{n}}.
\end{equation}

Then, we can bound the potential term. Using Equation~\eqref{eq:Y_n}, and $X_{n+1} = Y_{n+1} \circ (I - \gamma \nabla F)$, we have 
\begin{equation}
    \label{eq:explicitimplicit}
X_{n+1} = I - \gamma \nabla F - \gamma \nabla_W \cH(\mu_{n+1})(X_{n+1}).
\end{equation}
Since $F$ is $L$-smooth, we have~\cite{bau-com-livre11},
\begin{equation}
F(z) \leq F(x) + \ps{\nabla F(x),z-x} + \frac{L}{2}\|x-z\|^2, \qquad \forall\,\,x,z \in\sX.
\end{equation}
Replacing $z$ by $X_{n+1}(x)$, we obtain 
\begin{align}
F(X_{n+1}(x)) \leq F(x) - &\gamma \ps{\nabla F(x), \nabla F(x) + \nabla_W \cH(\mu_{n+1})( X_{n+1}(x))} \\
+&\frac{L\gamma^2}{2}\| \nabla F(x) + \nabla_W \cH(\mu_{n+1})(X_{n+1}(x))\|^2.
\end{align}
Integrating w.r.t. $\mu_n$, 
\begin{multline}
\label{eq:descent-F}
\cE_F(\mu_{n+1}) \leq \cE_F(\mu_{n}) -\gamma \ps{\nabla F, \nabla F + \nabla_W \cH(\mu_{n+1})(X_{n+1})}_{\mu_n} \\+ \frac{L\gamma^2}{2}\| \nabla F + \nabla_W \cH(\mu_{n+1})(X_{n+1})\|_{\mu_n}^2.
\end{multline}
Then, recalling~\eqref{eq:explicitimplicit} and summing equations ~\eqref{eq:descent-H} and ~\eqref{eq:descent-F}, we get
\begin{align*}
\cH(\mu_{n+1})+\cE_F(\mu_{n+1}) \leq& \cH(\mu_{n})+\cE_F(\mu_{n}) \\
&- \gamma \ps{\nabla_W \cH(\mu_{n+1})(X_{n+1}), \nabla F + \nabla_W \cH(\mu_{n+1})(X_{n+1})}_{\mu_{n}}\\
&- \gamma \ps{\nabla F, \nabla F + \nabla_W \cH(\mu_{n+1})(X_{n+1})}_{\mu_n} + \frac{L\gamma^2}{2}\| \nabla F + \nabla_W \cH(\mu_{n+1})(X_{n+1})\|_{\mu_n}^2\\
\leq& \cH(\mu_{n})+\cE_F(\mu_{n}) -\gamma\left(1- \frac{L\gamma}{2}\right)\| \nabla F + \nabla_W \cH(\mu_{n+1})(X_{n+1})\|_{\mu_n}^2.
\end{align*}

\section{Proof of Lemma~\ref{lem:evi-prox}}
Recall~\eqref{eq:Y_n},
\begin{equation*}
Y_{n+1} = I - \gamma \nabla_W \cH(\mu_{n+1})(Y_{n+1}) = T_{\nu_{n+1}}^{\mu_{n+1}}.
\end{equation*}
Since $(Y_{n+1}, T_{\nu_{n+1}}^{\pi})_{\#} \nu_{n+1}$ is a coupling between $\mu_{n+1}$ and $\pi$, we can upper bound the Wasserstein distance between $\mu_{n+1}$ and $\pi$ as:
\begin{align}\label{eq:evi-prox}
W^2(&\mu_{n+1},\pi) \nonumber
\leq  \|Y_{n+1} - T_{\nu_{n+1}}^{\pi}\|_{\nu_{n+1}}^2 \nonumber\\
= & \|I - T_{\nu_{n+1}}^{\pi}\|_{\nu_{n+1}}^2 -2\gamma \ps{\nabla_W \cH(\mu_{n+1})(Y_{n+1}), I - T_{\nu_{n+1}}^{\pi}}_{\nu_{n+1}} + \gamma^2 \|\nabla_W \cH(\mu_{n+1})(Y_{n+1})\|_{\nu_{n+1}}^2 \nonumber\\
= & \|I - T_{\nu_{n+1}}^{\pi}\|_{\nu_{n+1}}^2 -2\gamma \ps{\nabla_W \cH(\mu_{n+1})(Y_{n+1}), Y_{n+1} - T_{\nu_{n+1}}^{\pi}}_{\nu_{n+1}} - \gamma^2 \|\nabla_W \cH(\mu_{n+1})(Y_{n+1})\|_{\nu_{n+1}}^2.
\end{align}
where $\|I - T_{\nu_{n+1}}^{\pi}\|_{\nu_{n+1}}^2=W^2(\nu_{n+1},\pi)$. Moreover, using Lemma~\ref{lem:general-geo-cvx} with $\mu = \mu_{n+1}$ and $\nu = \nu_{n+1}$, $$-2\gamma \ps{\nabla_W \cH(\mu_{n+1})(Y_{n+1}), Y_{n+1} - T_{\nu_{n+1}}^{\pi}}_{\nu_{n+1}} \leq -2\gamma\left(\cH(\mu_{n+1}) - \cH(\pi)\right).$$ Plugging the latter inequality into~\eqref{eq:evi-prox}, we get the result.

\section{Proof of Lemma~\ref{lem:evi-grad}}
Since $(I - \gamma \nabla F, T_{\mu_n}^{\pi})_{\#} \mu_{n}$ is a coupling between $\nu_{n+1}$ and $\pi$, we can upper bound the Wasserstein distance between $\nu_{n+1}$ and $\pi$ as:
\begin{align}
\label{eq:evi-grad}
W^2(\nu_{n+1},\pi) 
&\leq  \|(I - \gamma \nabla F) - T_{\mu_{n}}^{\pi}\|_{\mu_{n}}^2 \nonumber\\
&=  \|I - T_{\mu_{n}}^{\pi}\|_{\mu_{n}}^2 -2\gamma \ps{\nabla F, I - T_{\mu_{n}}^{\pi}}_{\mu_n} + \gamma^2 \|\nabla F\|_{\mu_{n}}^2.
\end{align}
where $\|I - T_{\mu_{n}}^{\pi}\|_{\mu_{n}}^2=W^2(\mu_{n},\pi)$. Moreover, since $F$ is $\lambda$-strongly convex, we have:
\begin{equation}
F(x) \leq F(z) + \ps{\nabla F(x),x-z} - \frac{\lambda}{2}\|x-z\|^2, \qquad \forall \,\,x,z\in \sX.
\end{equation}
Replacing $z$ by $T_{\mu_{n}}^{\pi}(x)$ and multiplying by $2\gamma$, we obtain 
$$-2\gamma \ps{\nabla F(x), x - T_{\mu_{n}}^{\pi}(x)} \leq -2\gamma\left(F(x)- F \circ T_{\mu_{n}}^{\pi}(x)\right) - \gamma\lambda \|x - T_{\mu_{n}}^{\pi}(x)\|^2.$$
Integrating w.r.t. $\mu_n$ results in
$$-2\gamma \ps{\nabla F, I - T_{\mu_{n}}^{\pi}}_{\mu_n} \leq -2\gamma\left(\cE_F(\mu_n) - \cE_F(\pi)\right) - \gamma\lambda W^2(\mu_{n},\pi).$$
Plugging the latter inequality into~\eqref{eq:evi-grad} gives the result.

\section{Proof of Proposition~\ref{prop:evi}}
Recall that $Y_{n+1} = T_{\nu_{n+1}}^{\mu_{n+1}}$.
Combining Lemma~\ref{lem:evi-prox} and Lemma~\ref{lem:evi-grad}, we firstly get
\begin{align}
\label{eq:evi-not-finish}
W^2(\mu_{n+1},\pi) \leq& (1-\gamma\lambda)W^2(\mu_{n},\pi) \nonumber-2\gamma \left(\cE_F(\mu_n) + \cH(\mu_{n+1}) - \cE_F(\pi) - \cH(\pi)\right) \nonumber\\
&\quad + \gamma^2\|\nabla F\|_{\mu_n}^2 - \gamma^2\|\nabla_W \cH(\mu_{n+1})(Y_{n+1})\|_{\nu_{n+1}}^2.
\end{align}
Multiplying~\eqref{eq:descent-F} by $2\gamma$,
\begin{align*}
-2\gamma \cE_F(\mu_{n}) \leq& -2\gamma \cE_F(\mu_{n+1}) \\
&-2\gamma^2 \ps{\nabla F,\nabla_W \cH(\mu_{n+1})(X_{n+1})}_{\mu_n} - 2\gamma^2 \|\nabla F\|_{\mu_n} + L\gamma^3\| \nabla F + \nabla_W \cH(\mu_{n+1})(X_{n+1})\|_{\mu_n}^2.
\end{align*}
Moreover, using the transfer lemma, $\|\nabla_W \cH(\mu_{n+1})(X_{n+1})\|_{\mu_{n}}^2 = \|\nabla_W \cH(\mu_{n+1})(Y_{n+1})\|_{\nu_{n+1}}^2$. Therefore, 
\begin{align*}
-2\gamma& \cE_F(\mu_{n}) + \gamma^2\|\nabla F\|_{\mu_n}^2 - \gamma^2\|\nabla_W \cH(\mu_{n+1})(Y_{n+1})\|_{\nu_{n+1}}^2\\
\leq & -2\gamma \cE_F(\mu_{n+1}) - \gamma^2\|\nabla F\|_{\mu_n}^2 - \gamma^2\|\nabla_W \cH(\mu_{n+1})(X_{n+1})\|_{\mu_{n}}^2 -2\gamma^2 \ps{\nabla F,\nabla_W \cH(\mu_{n+1})(X_{n+1})}_{\mu_n}\\
& + L\gamma^3\| \nabla F + \nabla_W \cH(\mu_{n+1})(X_{n+1})\|_{\mu_n}^2\\
\leq & -2\gamma \cE_F(\mu_{n+1}) - \gamma^2(1-L\gamma)\|\nabla F + \nabla_W \cH(\mu_{n+1})(X_{n+1})\|_{\mu_n}^2.
\end{align*}
Plugging the last inequality into~\eqref{eq:evi-not-finish},
\begin{align*}
W^2(\mu_{n+1},\pi) \leq& (1-\gamma\lambda)W^2(\mu_{n},\pi) \nonumber-2\gamma \left(\cE_F(\mu_{n+1}) + \cH(\mu_{n+1}) - \cE_F(\pi) - \cH(\pi)\right) \nonumber\\
&- \gamma^2(1-\gamma L)\| \nabla F + \nabla_W \cH(\mu_{n+1})(X_{n+1})\|_{\mu_n}^2.
\end{align*}


\section{Proof of Corollary~\ref{cor:cvx}}
Let $\cL_n \eqdef 2\gamma n(\cG(\mu_n)-\cG(\mu_\star)) + W^2(\mu_n,\mu_\star)$.
From 
 Theorem~\ref{th:descent}, $\cG(\mu_{n+1})-\cG(\mu_\star) \leq \cG(\mu_n)-\cG(\mu_\star)$ if $\gamma < 1/L$. Therefore,
$$
2\gamma n (\cG(\mu_{n+1})-\cG(\mu_\star)) + 2\gamma(\cG(\mu_{n+1})-\cG(\mu_\star)) + W^2(\mu_{n+1},\mu_\star) \leq 2\gamma n (\cG(\mu_{n})-\cG(\mu_\star))+ W^2(\mu_n,\mu_\star),
$$
where we used Theorem~\ref{th:evi} with $\lambda = 0$ (recall that $\lambda \geq 0$). In other words,
$\cL_{n+1} \leq \cL_n$. Finally, 
$$
2\gamma n (\cG(\mu_{n})-\cG(\mu_\star)) \leq \cL_{n} \leq \cL_0 = W^2(\mu_0,\mu_\star).
$$

\section{Proof of Corollary~\ref{cor:scvx}}
Since the $\cG(\mu_{n+1})-\cG(\mu_\star)$ is nonnegative, from Theorem~\ref{th:evi},  
$$
W^2(\mu_{n+1},\mu_\star) \leq (1-\gamma\lambda)W^2(\mu_{n},\mu_\star).
$$
We get the result by iterating.


\end{document}